\documentclass[11pt]{amsart}
\usepackage{amsmath,mathrsfs}
\usepackage{amssymb}
\usepackage{latexsym}
\usepackage{amsfonts}
\usepackage{graphicx}
\usepackage{psfrag}

\newtheorem{Pa}{Paper}[section]
\newtheorem{theorem}[Pa]{{\bf Theorem}}

\newtheorem{corollary}[Pa]{{\bf Corollary}}

\newtheorem{proposition}[Pa]{{\bf Proposition}}
\newtheorem{observation}[Pa]{{\bf Observation}}

\newtheorem{example}[Pa]{{\bf Example}}

\def\C{\mathbb C}
\def\R{{\mathbb R}}

\begin{document}
\bibliographystyle{plain}
\title[Rational generalized positive functions and interpolation]
{Convex Cones of
Generalized Positive Rational Functions
and the Nevanlinna-Pick Interpolation}
\author[D. Alpay]{Daniel Alpay}
\address{(DA) Department of mathematics,
Ben-Gurion University of the Negev, P.O. Box 653, Beer-Sheva
84105, Israel} \email{dany@math.bgu.ac.il}
\author[I. Lewkowicz]{Izchak ~Lewkowicz}
\address{(IL) Department of electrical
engineering, Ben-Gurion University of the Negev, P.O.Box 653,
Beer-Sheva 84105, Israel} \email{izchak@ee.bgu.ac.il}

\thanks{This research is partially supported by BSF grant no
2010117}

\thanks{D. Alpay thanks the
Earl Katz family for endowing the chair which supported his
research. This research is part of the European Science
Foundation Networking Program HCAA}
\date{}



\maketitle
\mbox{Dedicated to appreciated colleagues Avraham Berman, Moshe
Goldberg and Raphael Loewy}
\begin{abstract}
Scalar rational functions with a non-negative real part on the
right half plane, called positive, are classical in the study of
electrical networks, dissipative systems, Nevanlinna-Pick
interpolation and other areas. We here study generalized positive
functions, i.e with a non-negative real part on the imaginary
axis. These functions form a Convex Invertible Cone, {\bf cic} in
short, and we explore two partitionings of this set: (i) into
(infinitely many non-invertible) convex cones of functions with
prescribed poles and zeroes in the right half plane and (ii)
each generalized positive function can be written as a sum
of even and odd parts. The sets of even generalized positive 
and odd 
%
%
functions form sub{\bf cic}s.

It is well known that the Nevanlinna-Pick interpolation problem 
is not always solvable by positive functions. Unfortunately,
there is no computationally simple procedure to carry out this
interpolation in the framework of generalized positive functions.
Through examples it is illustrated how the two above partitionings
of generalized positive functions can be exploited to introduce
simple ways to carry out the Nevanlinna-Pick interpolation.

Finally we show that only some of these properties are carried
over to rational generalized bounded functions, mapping the
imaginary axis to the unit disk.
\end{abstract}

\section{Introduction}
\subsection{Historical perspective}

Functions which are analytic in the open right half-plane
$\C_+$ and with a non-negative real part there
\begin{equation}\label{eq:pdef}
\begin{matrix}{\rm Re}~p(s)\ge 0,&~&s\in\C_+,\end{matrix}
\end{equation}
here denoted by ${\mathcal P}$, play an important in the theory
of electrical networks they were first studied around 1930 by W.
Cauer, \cite{Ca1}, \cite{Ca2} and O. Brune (who first coined the
name {\sl positive} for such functions; see \cite[Definition 1,
p.  25]{brune}, see also \cite{Br}). These functions also serve
as the corner stone of the theory of linear dissipative systems
(a.k.a absolutely stable), see e.g. \cite[Theorem 2.7.1]{AV},
\cite[3.18]{Be},
\cite{P} and \cite{wohlers}.

One can weaken condition \eqref{eq:pdef} and assume that a
function $~p~$ is analytic almost everywhere on the imaginary
axis satisfying
\[
\begin{matrix}{\rm Re}~p(s)\ge 0,&~&s\in{i}\R\end{matrix}.
\]
These functions will be
called {\sl generalized positive} and thus denoted by
$~\mathcal{GP}$. They were first addressed more than forty
years ago by B.D.O. Anderson and J.B. Moore in \cite{AM}.
\vskip 0.2cm

We shall denote by $~\C_-~$ the open left half of the
complex plane. We also denote the closed right half plane
by $~\overline{\C_+}~$ $(=i\R\cup\C_+)$.  We here consider
the field of scalar rational functions of a complex variable
$~s~$ with complex coefficients.  Throughout this work we
denote by $~\mathcal{GP}~$ the set of scalar ~{\em rational}~
generalized positive functions and by $~\mathcal{P}~$ its
subset of positive functions.
\vskip 0.2cm

In the sequel, we shall find it convenient to denote for an
arbitrary rational function $~g(s)$,
\[
g^{\#}(s) :=(g(-s^*))^*.
\]

The following result appeared in \cite{AL1}
(to help the reader we shall use $p\in{\mathcal P}$
and $\psi\in\mathcal{GP}$):

\begin{theorem}\label{Th:gp}
A rational function $~\psi(s)~$ is generalized positive if
and only if it admits the factorization
\begin{equation}\label{char}
\begin{matrix}\psi(s)=g(s)p(s)g^{\#}(s),&~&s\in\C,
\end{matrix}
\end{equation}
where $~p\in{\mathcal P}$ and $g$ is rational such that
both $~g~$ and $~g^{-1}~$ analytic in $\C_-~$ and
\begin{footnote}{Recall, the degree of a rational function is
taken to be the maximum between the degrees of the numerator
and the denominator.}\end{footnote}
\mbox{${\rm deg}(g)\in[0,~{\rm deg}(\psi)]$.}

Moreover, one can always\begin{footnote} {assuming
$\psi(s)\not\equiv 0$}\end{footnote} find $~s_o\in{i}\R$
so that in \eqref{char} the functions $~p(s)~$ and
$~g(s)~$ are uniquely determined by
$~0\not=\psi(s_o)=p(s_o)~$ and $~g(s_o)=1.$
\end{theorem}

Factorization results of the form of \eqref{char} are
well known in other frameworks. To name but three:

(i) Schur functions, analytically mapping the
open unit disk to its closure, rational generalized Schur
functions mapping the unit circle to the
closed unit disk, see e.g. \cite{Ta}, \cite{Ch}, \cite{Ak}.
Factorization result of generalized Schur
functions appeared in 
\cite[Theorem 3.2]{KL1}.

(ii) Carath\'{e}odory functions, mapping the open
unit disk to $~\overline{\C_+}$ and the rational
generalized (=pseudo) Carath\'{e}odory functions,
mapping the unit circle to $~\overline{\C_+}$.
Factorization result of generalized
Carath\'{e}odory functions appeared in
\cite[Theorem 3.1]{DGK}.

(iii) Nevanlinna functions, analytically mapping the closed
upper half plane to its closure and the rational generalized
Nevanlinna functions mapping the real axis to the closed upper
half plane. Factorization of generalized Nevanlinna functions
appeared in parallel in \cite{DHdS} and \cite{DLLS1} and further
explored and extended to operator valued functions in \cite{L1},
\cite{L2} and \cite{L3}.\\

An extended version of $\mathcal{GP}$ functions was introduced
by M.G. Kre\u{\i}n and H. Langer in a long and
celebrated series of papers of which we cite only \cite{KL1},
\cite{KL2} (note that they studied functions meromorphic in the
open upper half plane, or in the open disk).
\vskip 0.2cm

Important part of the existing research on $~\mathcal{GP}~$
functions is neither confined to scalar functions nor to rational
functions. Restricting the discussion to scalar rational functions
mapping the imaginary axis into the right half-plane, enabled us in
\cite{AL1} to obtain an elementary proof for the factorization
\eqref{char}. On the expense of generality, to keep the exposition
simple, we here adhere to the case of scalar rational functions.

\subsection{The current work}
Traditionally, $\mathcal{GP}$ functions were studied almost
uniquely by mathematicians. They were addressed in the framework of
not necessarily rational functions. In contrast, Electrical
engineers have long been interested in rational positive functions
(impedance of ~R-L-C~ electrical circuits).
\vskip 0.2cm

In this work we take a challenging task try to simultaneously
address both audiences. Thus, depending on their background, some
readers may find part
of the statements nearly obvious and others not clear at all.
Moreover, we try to address skeptical questions like:

(An engineer): "Why should I be interested in the extension
of positive functions to the $\mathcal{GP}$ framework?"

(A mathematician): "Why should I be interested in scalar
rational $\mathcal{GP}$ functions if the operator-valued
non-rational case has already been addressed ?"
\vskip 0.2cm

As already pointed out above, positive functions have been a
corner-stone in system theory. We thus believe that a sufficiently
good motivation for a researcher in this field to explore
$\mathcal{GP}$
functions, is to gain a proper perspective on the subset of
$\mathcal{P}$ functions. In particular to understand which of
the properties of $\mathcal{P}$ exist in the larger set
$\mathcal{GP}$ and which are peculiar to positive functions.
\vskip 0.2cm

As a prime example we point out that rational $\mathcal{GP}$
functions (bounded at infinity) may be characterized through the
Generalized Positive Real Lemma (in the positive case a.k.a. the
Kalman-Yakubovich-Popov Lemma) see e.g. \cite{AM}, \cite{DDGK}
for early accounts and most recently \cite{AL2}.
\vskip 0.2cm

One can examine this gap between  $\mathcal{GP}$ functions
and its Hurwitz stable subset of $\mathcal{P}$ functions,
from the Matrix Lyapunov Equation point of view. It has been of
interest to identify which of the
properties of the Lyapunov Equation for Hurwitz stable matrices
are carried over to the general inertia case, see e.g.
\cite[Chapter 2 and Section 4.4]{HJ2} or the relation to the
dimension of the controllable subspace in \cite{Lo}. Resorting
to the matrix Lyapunov equation was not just a metaphor, it is
part of the above mentioned Generalized Positive Real Lemma,
see e.g. \cite{AL2}, \cite{AM}, \cite{DDGK}.
\vskip 0.2cm

Recall that convex sets are an essential ingredient in
optimization. As an illustration, a typical control engineering
problem would be: ``Find among all stabilizing controllers the
one which minimizes a certain property". Such problems are
easier to address if the set of all controllers or closed
loop systems is convex. Convex sets also serve as a model for
uncertainty, e.g. the celebrated Kharitonov Theorem for 
checking the Hurwitz stability of a polytope of polynomials,
see \cite{Bar}. Hence, one is motivated in studying convex
sets of rational functions which are Hurwitz stable (and then
look for stable minimum phase). In Section \ref{sec:maximal}
we identify maximal convex cones of rational functions with
prescribed poles and zeroes in the open right (or left) half
plane.
\vskip 0.2cm

Using the above result, we introduce in Section \ref{sec:GP_g}
a partitioning of all $\mathcal{GP}$ functions into ``small",
yet maximal, convex cones denoted by $\mathcal{GP}_g$, with
prescribed poles and zeroes in the right half plane. Each of
these sets is a replica of $\mathcal{P}$ functions.
\vskip 0.2cm

Recall that a convex cone which in addition is closed under
inversion is called a Convex Invertible Cone, {\bf cic}~ in
short\begin{footnote}{Strictly speaking, this means that
whenever the inverse exists, it also belongs to the set, e.g.
the set of positive semidefinite matrices is a {\bf cic}. In
contrast, the open upper half of $\C$ is not.}
\end{footnote}, see e.g. \cite{CL1}, \cite{CL2} and \cite{CL4}.
It is easy to see that the set of $\mathcal{GP}$ functions is
closed under positive scaling, summation and inversion, i.e.
a {\bf cic} and ${\mathcal P}$ is a sub{\bf cic} of it. More
precisely, ${\mathcal P}$ is a maximal {\bf cic}~ of functions
which are analytic in $~\C_+$, see \cite[Proposition 4.1.1]{CL4}
and Proposition \ref{MaxG} below.
\vskip 0.2cm

In Section \ref{sec:odd} partition each $\mathcal{GP}$
function to a sum of {\em even}~ and {\em odd}
generalized positive
functions. It turns out that
the sets of {\em even}~ and {\em odd}~ generalized positive
functions are sub{\bf cic}s of $\mathcal{GP}$.
First, the
set of all $~{\mathcal Odd}~$ functions (i.e. generalized
lossless) which is of particular interest, is then studied.
In Section
\ref{sec:even} even $\mathcal{GP}$ functions are explored.
As a by-product of this partitioning it is shown that
it is only within the larger $\mathcal{GP}$ set that a positive
function can always be written as a sum of even and odd part.
\vskip 0.2cm

From an applications point of view, Nevanlinna-Pick type
interpolation was originally motivated by the design of the
driving point impedance of ~R-L-C~ electrical circuits and
subsequently by $~H_{\infty}$ control, both restricted to
positive functions. There are good reasons to study
Nevanlinna-Pick interpolation over $\mathcal{GP}$ functions:
\vskip 0.2cm

Recall that in the positive case solution exists, if and
only if, the Pick matrix associated with the data points is
positive semidefinite. In the $\mathcal{GP}$ framework, this
restriction is removed, namely for almost any set of data
points a solution exists\begin{footnote}{For generalized Schur
function this follows from \cite{DD} and
by appropriate Cayley transforms (of the functions and of the
variable) this is true for $\mathcal{GP}$ functions as well.
}\end{footnote}.
Note that in some applications, Hurwitz stability is not
required e.g. when the data is not associated with an
input-output dynamical system or when the system at hand is not
passive.
\vskip 0.2cm

Nevanlinna-Pick interpolation problem of generalized Schur and
Nevanlinna functions has been well addressed in the literature:
For generalized Nevanlinna functions see e.g. \cite{ABD},
\cite{AD}, \cite{ADLS}, \cite{bol_oam}, \cite[Section 3]{DZ} and
\cite{ADLRS}. For generalized Schur functions see e.g. \cite{AADLW},
\cite{Ball}, \cite{BH}, \cite{bol_oam2}, \cite{DD} and \cite{LW}
Nonetheless from computational point of view the known procedure
is involved.
\vskip 0.2cm

In each of the Sections: \ref{sec:GP_g}, \ref{sec:odd} and
\ref{sec:even}, we illustrate through examples how can one
exploit the new structural results to simplify the
Nevanlinna-Pick interpolation problem. A careful examination of
this examples suggests directions for future research. Some of
them are stated in Section \ref{sec:future}.
\vskip 0.2cm

Recall that a function is called ~{\em bounded},
denoted by $~f_b\in{\mathcal B}$, (commonly the real case is
addressed) if it analytically maps $~\C_+~$ to the closed unit
disk, see e.g. \cite[Chapter 7]{AV}, \cite[Section 6.5]{Be}
and $~f_{gb}\in\mathcal{GB}~$ is {\em generalized bounded}~
if it maps the imaginary axis to the closed
unit disk, see e.g. \cite{DDGK}. It is known that through the
Cayley transform positive functions may be identified with
bounded functions. In Section \ref{sec:bounded} properties of
generalized bounded functions, which do not trivially follow
from this Cayley transform are explored. In particular it is
shown that one cannot easily mimic Section \ref{sec:GP_g} to
obtain a partitioning of all rational generalized bounded
functions to a union of sets with prescribed poles and zeroes
outside the unit disk.

%
\section{Maximal convex sets of rational functions
with prescribed poles and zeroes in $~\C_+~$ or in $~\C_-$}
\label{sec:maximal}
\setcounter{equation}{0}

In this section we consider poles and zeroes of sums of rational
functions. Up to possible cancellations, poles of a sum are
the union of the poles original functions. However, in general
little can be said about the zeroes of a sum. We now characterize
maximal convex sets of rational functions with prescribed poles
and zeroes in $~\C_+~$ (or in $~\C_-$).

We begin with some preliminaries. We find it convenient to
define the following sets
\begin{equation}\label{eq:G_+}
\begin{matrix}
{\mathcal G}_-:=~{\rm all~rational~functions~with~poles~
and~zeroes~in}~\C_-~,\\
{\mathcal G}_+:=~{\rm all~rational~functions~with~poles~
and~zeroes~in}~\C_+~.
\end{matrix}
\end{equation}
Note that $g\in{\mathcal G}_-~$ is equivalent to
$~g^{\#}\in{\mathcal G}_+$.

\begin{example}\label{Ex:G_+^1}
{\rm
All degree one real functions with poles and zeroes in the
$\overline{\mathbb{C}_+}$
are given by
\begin{equation}\label{eq:G_a}
\mathcal{G}_1=\left\{ \frac{as-b}{cs-d}~:~ab\geq 0,~cd\geq 0,~
ad\not=bc.~\right\}.
\end{equation}
Now $~\mathcal{G}_1^{\#}=\left\{ \frac{as+b}{cs+d}~:~
ab\geq 0,~cd\geq 0,~ad\not=bc.~\right\}$,~ and all degree one
real functions within $\mathcal{P}$ is a subset of
$~\mathcal{G}_1^{\#}~$ where in addition $~a,~c\geq 0$.
Indeed, up to inversion, all degree one real functions in
$\mathcal{P}$ are of the form $~as+b~$ or $~\frac{a}{s}+b~$
with $a,b\geq 0$ (in electrical circuits terminology, the
driving point impedance of ~R-L~ or ~R-C~ networks,
respectively).
}\qed
\end{example}
\vskip 0.2cm

Next, for a {\em given} $~g_+~$ in $~\mathcal{G}_+$, see
\eqref{eq:G_+}, let $~\tilde{\mathcal G}_{g_+}~$ be the set of
all rational functions with prescribed poles and zeroes in
$~\C_+~$ given by,
\begin{equation}\label{TildeG}
\tilde{\mathcal G}_{g_+}:=\left\{\frac{n_-(s)}{d_-(s)}
g_+(s)~:~n_-(s), d_-(s)~ {\rm polynomials~with~roots~in}~
\overline{\C_-}\right\}.
\end{equation}
Obviously,
if in \eqref{TildeG} $g\in\tilde{\mathcal G}_{g_+}~$ satisfies
$\frac{n_-(s)}{d_-(s)}\equiv{\rm const}.~$
then in fact $~g\in\mathcal{G}_+~$ \eqref{eq:G_+}.
The set $~\tilde{\mathcal G}_{g_+}~$ can not be convex
as both $~-g_+(s)~$ and $~sg_+(s)~$ belong to it, but their sum
has an additional zero in $~\C_+$. Yet, it is of interest to
identify maximal convex subsets of $~\tilde{\mathcal G}_{g_+}$
in \eqref{TildeG}. To this end, denote by $~{\mathcal G}_o~$
the set of rational functions, which along with their inverses,
are analytic in both open half planes. Namely, whenever
$~g_o\in{\mathcal G}_o$, it is of the form
\begin{equation}\label{g_o}
g_o(s)=c\prod\limits_{j=1}^m(s-ir_j)^{\eta_j}~,
\end{equation}
where $~r_j\in\R~$ are all distinct, $~\eta_j~$ are integers
(positive or negative) and $~c\in\C$. We shall use the convention
that $~\prod\limits_1^0=1$, so that also
$~g(s)\equiv {\rm const.}$ belongs to this set.
\vskip 0.2cm

For three given functions $~g_+\in\mathcal{G}_+~$,
$~g_-\in\mathcal{G}_-$ see \eqref{eq:G_+}, and $~g_o~$ see
\eqref{g_o}, define,
\begin{equation}\label{eq:G}
{\mathcal G}_{g_+,~g_-,~g_o}:=\{g_+(s)p(s)g_o(s)g_-(s)~:~
p\in{\mathcal P}~\}.
\end{equation}
By construction, in $~\C_+~$ the poles and zeroes of all
functions\begin{footnote}{excluding the zero function}
\end{footnote} in $~{\mathcal G}_{g_+,~g_-,~g_o}~$ are exactly
those
of $~g_+(s)$. On $~i\R~$ poles and zeroes are almost fixed
in the following sense. Considering \eqref{g_o}, functions in
$~{\mathcal G}_{g_+,~g_-,~g_o}~$ have factors of the form
$~(s-ir_j)^{\eta_j+1}$, $~(s-ir_j)^{\eta_j-1}~$ or
$~(s-ir_j)^{\eta_j}~$ depending on $~p(s)~$ having at $~s=ir_j$,
a zero, a pole, neither zero nor pole, respectively.
\vskip 0.2cm

%
%
%
%
%
We can now describe convex sets of functions with prescribed
poles and zeroes in $~\C_+~$ and on $~i\R$, almost prescribed
in the above sense.

\begin{proposition}\label{MaxG}
The following is true:

\begin{itemize}

\item[(i)~~~~]{}The set $\mathcal{P}$ is a maximal convex
invertible cone, {\bf cic}, of rational functions analytic in
$~\C_+$.

\item[(ii)~~~]{}${\mathcal G}_{g_+,g_-,g_o}~$ in \eqref{eq:G},
is a maximal convex set of rational functions with prescribed
poles and zeroes in $~\C_+$.\\
In fact, if $~\phi\not\in{\mathcal G}_{g_+, g_-, g_o}$,
one can always find $~\psi\in{\mathcal G}_{g_+, g_-, g_o}~$
so that $~(\phi+\psi)\not\in\tilde{\mathcal G}_{g_+}~$
\eqref{TildeG}.
\end{itemize}
\end{proposition}
\vskip 0.2cm

The fact that the set ${\mathcal P}$ is a convex invertible cone,
{\bf cic}, is classical, see e.g. \cite[5.6]{Be}. The result in
item (i) is a small variation of \cite[Proposition 4.1.1]{CL4}.
\vskip 0.2cm

\noindent
{\bf Proof :}\quad (i) We first show that if $~h(s)~$ is a non-positive
function, one can always find a positive function $~p~$ so that
$~h+p~$ has a zero in $~\C_+$. Indeed, let $~h\not\in{\mathcal P}~$
be given. By definition there are points in $~\C_+~$ which are
mapped by $~h(s)~$ to $~\C_-$, i.e. there exist $~\alpha,
\gamma>0$, $~\beta, \delta\in\R~$ so that
\mbox{$h(s)_{|_{s=\alpha+i\beta}}=- \gamma+i\delta$}. Take now
$~p(s)=\frac{\gamma}{\alpha}s-i(
\frac{\beta\gamma}{\alpha}+\delta)$. Then clearly $~p\in{\mathcal
P}~$ and \mbox{$(p+h)(s)_{ |_{s=\alpha+i\beta}}=0$}, i.e. a zero
in $~\C_+~$. Next, note that $~\frac{1}{p+h}~$ is not analytic in
$~\C_+$. Since $~{\mathcal P}~$ is closed under inversion (i.e.
$~p\in{\mathcal P}~$ is equivalent to $~\frac{1}{p}\in{\mathcal P}$),
this part is established.\\
(ii) Let $~g_-\in{\mathcal G}_-~$, $~g_+\in{\mathcal G}_+~$ see
\eqref{eq:G_+} and $~g_o\in{\mathcal G}_o$, see \eqref{g_o}, be
given and let $~\phi(s)~$ be a rational function not in
$~{\mathcal G}_{g_+, g_o, g_-}$. To avoid triviality assume that
\mbox{$\phi\in\{\tilde{\mathcal G}_{g_+}\smallsetminus{\mathcal
G}_{g_+, g_-, g_o}\}$.} Next, denote,
\mbox{$~h(s)=g_o(s)^{-1}g_+(s)^{-1}\phi(s)g_-(s)^{-1}$.} Then by
\eqref{eq:G}, $~h\not\in{\mathcal P}$ (else $~\phi~$ would
have been in $~{\mathcal G}_{g_+, g_-, g_o})$. Take now
$~\psi(s)=g_+(s)g_o(s)p(s)g_-(s)~$ with the above $~g_+(s)$,
$g_o(s)$, $~g_-(s)~$ and $~p(s)~$ as in part (a) of this
proof. By construction, $~p\in{\mathcal P}~$ and thus
\mbox{$\psi\in{\mathcal G}_{g_+, g_-, g_o}$,} but
\mbox{$(\phi+\psi )\not\in\tilde{\mathcal G}_{g_+}$} since
this function has an additional zero in $~\C_+$.
(If this additional zero coincides with an existing pole, both in
$~\C_+$, still $~ (\phi+\psi)\not\in\tilde{\mathcal G}_{g_+}$).
Thus, this part of the claim is established and the proof is
complete. \qed\mbox{}
\vskip 0.2cm

As an illustration we have the following.
\begin{example}
{\rm
Take $g_+(s)$, $g_o(s)$ and $g_-(s)$ in \eqref{eq:G} to be fixed
polynomials. A maximal (up to scaling) convex set of polynomials
whose roots in $\C_+$ are those of $g_+$, is given by
\[
\{ g_+(s)(s+a)g_o(s)g_-(s)~:~a\in\overline{{\mathbb C}_+}~\}.
\]
Indeed, $p(s)=s+a$, $a\in\overline{{\mathbb C}_+}$ are the only
polynomials in $\mathcal{P}$.
}
\qed
\end{example}

We conclude this section by stating the analogous results for
the left half plane. First, we denote by ${\mathcal P}^{\#}$ the
set of~ {\em para-positive} functions,
\begin{equation}\label{PP}
{\mathcal P}^{\#}:=\{~\psi~:~\psi^{\#}\in{\mathcal P}~\}.
\end{equation}
Thus, functions in ${\mathcal P}^{\#}$ map $~\overline{\C_-}~$
to $~\overline{\C_+}$. In particular,
$~{\mathcal P}^{\#}\subset\mathcal{GP}$.
We can now state results which are dual to Proposition \ref{MaxG}.

\begin{observation}\label{C_-}
The following is true:

\begin{itemize}
\item[(i)~~~~]{}The set ${\mathcal P}^{\#}$ is a maximal convex
invertible cone, {\bf cic}, of rational functions analytic in $~\C_-$.

\item[(ii)~~~]{}Let $~g_-\in{\mathcal G}_-$, $~g_+\in{\mathcal G}_+~$
and $~g_o\in{\mathcal G}_o~$ be given. The set
$$\{g_+(s)p^{\#}(s)g_o(s)g_-(s)~:~p\in{\mathcal P}~\}$$
is a maximal convex set whose poles and zeroes in $~\C_-~$
are precisely those of $~g_-(s)$.
\end{itemize}
\end{observation}

\section{Convex partitioning of $~\mathcal{GP}~$ functions}
\label{sec:GP_g}
\setcounter{equation}{0}

We now address ourselves to subsets of generalized positive
functions within $~{\mathcal G}_{g_+,~g_-,~g_o}~$ in \eqref{eq:G},
namely sets of the form
$~{\mathcal G}_{g_+,~g_-,~g_o}\cap\mathcal{GP}$. To this end,
we introduce the
following set, using \eqref{eq:G_+} and \eqref{g_o},
\begin{equation}\label{overline{G}_+}
\overline{{\mathcal G}_+}:=\{ g_+(s)g_o(s)~:~g_+\in{\mathcal G}_+,~
g_o\in{\mathcal G}_o~\}.
\end{equation}
Note that $~g\in\overline{{\mathcal G}_+}~$ means that both $~g~$
and $~g^{-1}~$ are analytic in $~\C_-$.
For example, all degree one functions in $\overline{{\mathcal G}_+}$
are given by (the real subset $\mathcal{G}_1$ was described in
\eqref{eq:G_a}),
\begin{equation}\label{eq:G_+^1}
\hat{\mathcal G}
=\left\{ \frac{as-b}{cs-d}~:~{\rm Re}(a^*b)\geq 0,~
{\rm Re}(c^*d)\geq 0,~ad\not=bc\right\}.
\end{equation}
Using this notation, we shall hereafter simply write
$$\mathcal{GP}_g:={\mathcal G}_{g_+,~g_-,~g_o}
\cap\mathcal{GP}.$$
By Theorem \ref{Th:gp}, for given $~g\in\overline{{\mathcal G}_+}$
this can be equivalently written as
\begin{equation}\label{eq:GP_g}
\mathcal{GP}_g=\{~gpg^{\#}~:~p\in{\mathcal P}\}.
\end{equation}
For a given $~g\in\overline{{\mathcal G}_+}$, the set
$~\mathcal{GP}_g$ is a replica of $~{\mathcal P}$. 
Nevertheless, the picture in $~\mathcal{GP}_g$ is richer.
\begin{example}\label{Ex:GP_g}
{\rm
We here illustrate two properties of the set $~\mathcal{GP}_g$
\eqref{eq:GP_g}, where $~g\in\overline{{\mathcal G}_+}~$ is
given:

(a) If $~\psi_j=gp_jg^{\#}$ with $~g~$ fixed,
$~{\rm deg}~p_1>{\rm deg}~p_2$,
does not always imply $~{\rm deg}~\psi_1>{\rm deg}~\psi_2~$.\\
(b) In this set, it is only in $~\C_+~$ that the poles and
zeroes are fixed. On $~i\R~$ they are almost prescribed,
and in $~\C_-~$ they are not fixed.

Take $~g(s)=\frac{s-2}{s}~$ ($g\in{\mathcal G}_1$, see
\eqref{eq:G_a}). Thus, \mbox{$\mathcal{GP}_g=\left\{~\frac{s^2-4}
{s^2}p(s)~:~p\in{\mathcal P}\right\}$}. We here mention, but
five interesting samples,
\begin{center}
$\begin{matrix}~&~&p(s)&~&~&\psi(s)=g(s)p(s)g^{\#}(s)\\ ~&~&~&~\\
{\rm (i)}&~&\frac{s}{s+2}&~&~&\frac{s-2}{s}\\ 
{\rm (ii)}&~&\frac{s}{(s+2)^2}&~&~&\frac{s-2}{s(s+2)}\\
{\rm (iii)}&~&1&~&~&\frac{s^2-4}{s^2}\\
{\rm (iv)}&~&\frac{s+2}{s}&~&~&\frac{(s+2)^2(s-2)}{s^3}\\
{\rm (v)}&~&\frac{s(s+2i)}{s+i}&~&~&\frac{(s^2-4)(s+2i)}{s(s+i)}~.
\end{matrix}$
\end{center}
(a) These five functions are ordered so that the degree of
$~\psi_j(s)~$ is non-decreasing. In contrast, the degree of
the corresponding $~p_j(s)~$ ``fluctuates".

(b) In $~\C_+$, there is always a zero with a unit multiplicity at
$~s=+2$.\\
On $~i\R$, there is a pole at the origin. Its generic multiplicity
is two, but it may also be one or three (i.e. at the origin $~p~$ has
a pole e.g. (iv), no pole nor zero e.g. (iii), or a zero e.g.
(i), (ii), (v) respectively). $~\psi~$ may have additional poles
or zeroes, see e.g. (v).\\
In $~\C_-~$ poles and zeroes are not fixed, see e.g. the point
$~s=-2$.
}\qed
\end{example}

Using the notation of \eqref{overline{G}_+} Theorem
\ref{Th:gp} may be formulated as saying that having
$~\psi\in\mathcal{GP}$ is equivalent to
\mbox{$\psi(s)=g(s)p(s)g^{\#}(s)$} for some
$~g\in\overline{{\mathcal G}_+}~$ and some $~p\in{\mathcal P}$.
Thus, we can use the last result to introduce a {\em convex}~
partitioning of all $~\mathcal{GP}~$ functions. The proof is
immediate and thus omitted.
\begin{observation}\label{ConvexPartitioning}
Let $~\mathcal{GP}_g~$ be as in \eqref{eq:GP_g}. Then,
\begin{itemize}
\item[(i)~~~~]{}
$\mathcal{GP}_g$ is a convex sub-cone of $\mathcal{GP}$.

\item[(ii)~~~]{}For $~g_1, g_2\in\overline{{\mathcal G}_+}~$
$~\mathcal{GP}_{g_2}=
\left(\frac{g_2}{g_1}\right)\mathcal{GP}_{g_1}
\left(\frac{g_2}{g_1}\right)^{\#}$.

\item[(iii)~~]{}Let $~g_1, g_2\in\overline{{\mathcal G}_+}~$ be
so that $~g_2\not\equiv cg_1$, for some constant $~c$, then,
\mbox{$\mathcal{GP}_{g_1}\cap\mathcal{GP}_{g_2}=\{ 0\}$.}

\item[(iv)~~]{}$\mathcal{GP}=\bigcup\limits_{
g\in\overline{{\mathcal G}_+}}\mathcal{GP}_g~$.

\item[(v)~~~]{}$(\mathcal{GP}_g)^{-1}=\mathcal{GP}_{(g^{\#})^{-1}}$.
\end{itemize}
\end{observation}

{\bf Proof}. Items (i), (ii) and (v) are immediate from
\eqref{eq:GP_g}. Item (iv) follows from Theorem \ref{Th:gp} 
along with \eqref{eq:GP_g}.

As to item (iii),
assume that there exits $~\psi~$ within
$\mathcal{GP}_{g_1}\cap\mathcal{GP}_{g_2}$ for some
$~g_1,g_2\in\overline{{\mathcal G}_+}$.
We shall find it
convenient to factorize $~g_j=g_{o,j}g_{+,j}$ with $j=1,2$.
where $~g_{o,1},g_{o,2}\in\mathcal{G}_o$,
see \eqref{g_o}, and $~g_{+,1},g_{+,2}\in\mathcal{G}_+$, see
\eqref{eq:G_+}. As poles and zeroes of $~\psi~$ in $~\C_+~$ 
are uniquely determined by $~g_+$, without loss of generality
one can write $~g_1=g_{o,1}g_+$ and $~g_2=g_{o,2}g_+$ for some
$~g_+\in\mathcal{G}_+$.

Next assume that for $j=1,2$ and some $r\in\R$, $~g_{o,j}(s)$
have factors $(s-ir)^{m_j}$ and the corresponding $~p_j(s)~$
have factors $(s-ir)^{l_j}$, where $~m_j, l_j$ are integers
(not necessarily positive). Then in $\psi_j(s)$ the respective
factors are
\[
(s-ir)^{m_j}(s-ir)^{l_j}\left((s-ir)^{m_j}\right)^{\#}=
(-1)^{m_j}(s-ir)^{2m_j+l_j}.
\]
This implies that: (i) $~2m_1+l_1=2m_2+l_2$ and (ii)
$~m_1-m_2=2k$ for some integer $k$. Namely, $l_2-l_1=2(m_1-m_2)=4k$.
Now recall that on the imaginary axis poles and zeroes of positive
functions are simple, see e.g. \cite[Theorem 2.2]{AL1}, i.e.
$1\geq |l_j|$. This implies that $2\geq |l_2-l_1|$. But, 
$m_1\not=m_2$ implies that $|l_2-l_1|\geq 4$. Hence,
one can now conclude that $m_1=m_2$,
$l_1=l_2$ and since $r$ was arbitrary (up to a constant)
$g_{0,1}=g_{0,2}$, which in turn implies $g_1=g_2$ (up to a constant)
and the proof is complete.
\qed\mbox{}
\vskip 0.2cm

This partitioning of $\mathcal{GP}$ functions rightfully seems 
straightforward. In contrast, at the end of Section
\ref{sec:bounded}, we show that partitioning of $\mathcal{GB}$,
generalized bounded functions (or generalized Schur functions) in
the spirit of Observation \ref{ConvexPartitioning}, can not be
easily mimicked.
\vskip 0.2cm

For a given $~g\in\overline{{\mathcal G}_+}~$ we now wish to
identify minimal degree functions within $~\mathcal{GP}_g~$.

\begin{proposition}\label{MinDeg}
The following is true.

\begin{itemize}
\item[(i)~~~~]{}$g\in\overline{{\mathcal G}_+}~$ can
always be factored as
\mbox{$g(s)=c\prod\limits_{j=1}^q\psi_j(s)$} with
$~\psi_j^{\#}(s)~$ positive, see \eqref{PP}, and $~c\in\C$.

\item[(ii)~~~]{}Among all possible factorizations of
$~g\in\overline{{\mathcal G}_+}~$ let us choose
$~g(s)=g_1(s)g_2(s)~$ so that $~g_1^{\#}(s)~$ is positive
and $~{\rm deg}(g_1)~$ is maximal.\\
Then, $~gg_2^{\#}~$ is the minimal degree function in
$\mathcal{GP}_g$.
\end{itemize}
\end{proposition}

{\bf Proof :}\quad Item (i) is immediate from item (i) in
Observation \ref{C_-}. Specifically,
$~g\in\overline{{\mathcal G}_+}~$ can always be written as
$~g=\frac{\prod\limits_j(s-z_j)}
{\prod\limits_k(s-\pi_k)}~$ with $~z_j, \pi_k\in\overline{\C_+}$.
Note now that $~(s-z_j)^{\#}=-(s+z_j^*)~$ and
$~\frac{1}{(s-\pi_k)^{\#}}=\frac{-1}{s+\pi_k^*}~$. Namely,
up to sign change, this $~g^{\#}~$ is a product of degree
one positive functions.\\
Item (ii) stems from item (i) and \eqref{eq:GP_g} noting that
$~\psi\in\mathcal{GP}_g~$ can always be written as,
\[
\psi(s)=g(s)p(s)g^{\#}(s)=
g_1(s)g_2(s)p(s)g_1^{\#}(s)g_2^{\#}(s),
\]
so that $~g_1^{\#}(s)~$ is positive, see e.g. Example \ref{Ex:G_+^1}.
Thus, to reduce the degree of the above $~\psi(s)~$ choose
$~p(s)=\frac{1}{g_1^{\#}(s)}~$ so that
\[
\psi(s)=
g_1(s)g_2(s)p(s)g_1^{\#}(s)g_2^{\#}(s)_{|_{p=\frac{1}{g_1^{\#}}}}
=g_1(s)g_2(s)g_2^{\#}(s)=g(s)g_2^{\#}(s)
\]
and the construction is complete.
\qed
\vskip 0.2cm

The above construction is illustrated in part (a) of Example
\ref {Ex:GP_g}.
\vskip 0.2cm

Within the set $~\mathcal{P}$, the Nevanlinna-Pick interpolation
problem is classical, see e.g. \cite[Chapter 18]{BGR}, \cite{YS}.
Within $~\mathcal{GP}$, variants of this interpolation problem
are well studied, see e.g. \cite{ABD}, \cite{AD}, \cite{ADLS},
\cite{bol_oam}, \cite[Section 3]{DZ} and \cite{ADLRS} for
generalized Nevanlinna functions and for generalized Schur
functions see e.g. \cite{Ball}, \cite{BH}, \cite{bol_oam2}, and
\cite{DD} and \cite{AADLW}. Nonetheless from computational point
of view the procedure is involved. As an intermediate step, in
the following example we illustrate the fact that within the set
$\mathcal{GP}_g$, namely when $g\in\overline{{\mathcal G}_+}$ is
fixed, the Nevannlin-Pick interpolation problem 
reduces to the classical version within the set $~\mathcal{P}$,
which is computationally well established.

\begin{example}\label{Ex:GenNevPickGP_g}
{\rm
We here illustrate a Nevanlinna-Pick interpolation scheme within
the set $~\mathcal{GP}_g~$.
We look for $~\psi\in\mathcal{GP}_g~$ so that
\[
\psi(s)_{|_{s=1}}=1
\quad\quad{\rm and}\quad\quad
\psi(s)_{|_{s=2}}=4.
\]
(As the associated Pick matrix is {\mbox{\tiny$\begin{pmatrix}
1&\frac{5}{3}\\ \frac{5}{3}&2\end{pmatrix}$}}, its
determinant is negative, so there is no $\psi\in\mathcal{P}$).
\vskip 0.2cm

Take\begin{footnote}{Recall that we assume $g$ is prescribed.}
\end{footnote} $~g(s)=\frac{4}{7-3s}~$ i.e. a right half plane
pole at $~s=\frac{7}{3}~$ and a zero at infinity. Thus,
\mbox{$g(s)g^{\#}(s)=\frac{16}{49-9s^2}~$.}~ First denote
\begin{center}
$w_1:=g(s)g^{\#}(s)_{|_{s=1}}=${\mbox{\tiny$\frac{2}{5}$}}
$\quad\quad{\rm and}\quad\quad w_2:=g(s)g^{\#}(s)_{|_{s=2}}=$
{\mbox{\tiny$\frac{16}{13}$}}~.
\end{center}
Exploiting the $~\mathcal{GP}_g~$ structure, see \eqref{eq:GP_g}, we
actually seek $p\in{\mathcal P}$ so that
\begin{center}
$p(s)_{|_{s=1}}=$
{\mbox{\tiny$\frac{1}{w_1}$}}
$=${\mbox{\tiny$\frac{5}{2}$}}
$\quad\quad{\rm and}\quad\quad p(s)_{|_{s=2}}=$
{\mbox{\tiny$\frac{4}{w_2}$}}
$=$
{\mbox{\tiny$\frac{13}{4}$}}~.
\end{center}
But this is a classical Nevanlinna-Pick problem
and the resulting Pick matrix is, $\Pi=${\mbox{\tiny$
\begin{pmatrix}\frac{5}{2}&\frac{23}{12}\\ 
\frac{23}{12}&\frac{13}{8} \end{pmatrix}$}}.  As this
$~\Pi~$ is positive definite, there are infinitely many
solutions (using $\Pi$, they can all be parameterized, see
e.g. \cite[Chapter 18]{BGR}).

Take for instance two degree three interpolating functions,
\mbox{$p_a(s)=\frac{3s(s^2+9)}{4(s^2+2)}$} and

\mbox{$p_b(s)=\frac{9}{4}\left(1+\frac{s^3}{3(s^2+2)}\right)$.}
The resulting interpolating functions in $~\mathcal{GP}_g~$ are
$$\begin{matrix}
\psi_a(s)&=&g(s)p_a(s)g(s)^{\#}&=&
\frac{12(s^2+9)}{(s^2+2)(49-9s^2)}~,\\
\psi_b(s)&=&g(s)p_b(s)g(s)^{\#}&=&
\frac{12(s^3+3s^2+6)}{(s^2+2)(49-9s^2)}~.
\end{matrix}$$
Both $~\psi_a(s)~$ and $~\psi_b(s)~$ are of degree four.

Recall that the set of interpolating functions is convex.
Take for example
\begin{center}
$\psi_c(s):=${\mbox{\tiny$\frac{2}{9}$}}$\psi_a(s)+$
{\mbox{\tiny$\frac{7}{9}$}}$\psi_b(s)=g(s)
=\frac{4}{7-3s}~,$
\end{center}
to obtain a {\em unity~ degree}~ interpolating function within
$\mathcal{GP}_g$. From  item (ii) in Proposition
\ref{MinDeg} it follows that in fact this is the minimal degree
function within $\mathcal{GP}_g$ and in particular the minimal
degree interpolating function.
}\qed
\end{example}

\section{Odd functions - a sub{\bf cic} of generalized positive
functions}\label{sec:odd}
\setcounter{equation}{0}

As already mentioned, it is easy to see that the set of
$\mathcal{GP}$ functions is closed under positive scaling,
summation and inversion, i.e. a {\bf cic}. In the previous
section we introduced a partitioning of this {\bf cic} into
(infinitely many non invertible) convex cones of the form
$\mathcal{GP}_g$. We now explore a partitioning of each
generalized positive function into even and odd parts. It
turns out that the sets of even and odd generalized positive
functions are two sub{\bf cic}s of $\mathcal{GP}$.
\vskip 0.2cm

Abusing the real case terminology, for a given rational
function $f(s)$ we shall define the {\em even}~ and
{\em odd}~ parts as
\begin{equation}
\begin{matrix}\label{evenoddDef}
f_{\rm even}(s):=\frac{1}{2}\left(f(s)+f^{\#}(s)\right)&~&
f_{\rm odd}(s):=\frac{1}{2}(f(s)-f^{\#}(s)).
\end{matrix}
\end{equation}
Then, we also define the sets of all even and all odd
functions,
\begin{equation}\label{EvenOddDef}
\begin{matrix}
{\mathcal Even}:=\{~f(s)~:~f=f_{\rm even}~\}&~&
{\mathcal Odd}:=\{~f(s)~:~f=f_{\rm odd}~\}.
\end{matrix}
\end{equation}
The following observations are almost obvious, they are stated
for a comparison in the sequel.
\begin{proposition}\label{EvenOdd}
Let the sets $~{\mathcal Even}~$ and $~{\mathcal Odd}~$ be as in
\eqref{evenoddDef} and \eqref{EvenOddDef}.
\begin{itemize}
\item[(i)~~~~]{}${\mathcal Even}~$ and $~{\mathcal Odd}~$ are
convex invertible cones, {\rm{\bf cic}s} of rational functions.

\item[(ii)~~~]{}Let $~f,g~$ be rational functions.
If $~g\in{\mathcal Even}~$ then,
$$\begin{matrix}
(fg)_{\rm even}=f_{\rm even}g&~&(fg)_{\rm odd}=f_{\rm odd}g.
\end{matrix}$$
Conversely, if $~(fg)_{\rm even}=f_{\rm even}g~$ and
$~f\not\equiv 0$,
then~ \mbox{$~g\in{\mathcal Even}$.}
\end{itemize}
\end{proposition}

\noindent
{\bf Proof~}\quad
(i) Indeed, positive scaling and summation
are obvious. As to inversion note that if
$~f(s)=\frac{n(s)}{d(s)}$, with $~n,~d~$ polynomials, then
$~f\in{\mathcal Even}~$ is equivalent to $~nd^{\#}=n^{\#}d$,
which in turn means, $~f^{-1}\in{\mathcal Even}$. The reasoning
for  $~f\in{\mathcal Odd}~$ is similar and thus omitted.\\
(ii) A straightforward calculation shows that having
$~(fg)_{\rm even}=f_{\rm even}g~$ is equivalent
to $~f^{\#}g=f^{\#}g^{\#}$, which in turn for $~f\not\equiv 0~$
means $~g\in{\mathcal Even}$.\qed
\vskip 0.2cm

From item (i) in Proposition \ref{MaxG} and item (i) in Proposition
\ref{EvenOdd} it respectively follows that ${\mathcal P}$ and
${\mathcal Odd}$ are {\bf cic}s. Recall that a non-empty
intersection of {\bf cic}s is a {\bf cic}, see \cite[Observation
2.1]{CL4}. Thus, of particular interest is the sub{\bf cic} of all
positive odd rational functions
$~\mathcal{PO}:={\mathcal P}\cap{\mathcal Odd}$. The real subset
of $\mathcal{PO}$ functions are sometimes called ``lossless", or
``Foster" and they correspond to L-C circuits, see e.g. \cite{AV},
\cite{Be}. 
\vskip 0.2cm

$\mathcal{PO}$ functions can be parameterized, see e.g. \cite[5.13]{Be},
as
\begin{equation}\label{po}
\mathcal{PO}:=\left\{
p(s)=ir_o+a_os+\sum\limits_{j\geq 1}\frac{a_j}
{s-ir_j}~:~a_o\geq 0,~a_j>0,~r_j\in\R\right\}.
\end{equation}
Combining
Theorem \ref{Th:gp} together with \eqref{po} we have the following.

\begin{proposition}\label{gpo}
Let the set of odd functions, ${\mathcal Odd}$, be as in
\eqref{evenoddDef} and \eqref{EvenOddDef}. The following
statements are true.

\begin{itemize}
\item[(i)~~~~]{}${\mathcal Odd}=\mathcal{GP}\cap{\mathcal Odd}.$

\item[(ii)~~~]{} $\psi\in{\mathcal Odd}~$ if and only if $~\psi~$
maps the imaginary axis to itself.

\item[(iii)~~]{} The set $~{\mathcal Odd}~$ is closed (excluding the
zero function) under:\\
(a) real scaling, (b) addition, (c) inversion (d) composition and
(e) the product of an odd number of elements, i.e.
$~\prod\limits_{j}^{2m+1}\psi_j(s)~$ with $~\psi_j\in{\mathcal Odd}$,
$~m=0,~1,~\ldots$

\item[(iv)~~]{} $\psi\in{\mathcal Odd}~$ can always be written as
$$\psi(s)=g(s)\left( ir_o+a_os+\sum\limits_{j\geq 1}
\frac{a_j}{s-ir_j}\right)g^{\#}(s),$$
with $~a_o\geq 0,~a_j>0,~r_j\in\R~$ and
$~g\in{\overline{\mathcal G}_+}$ \eqref{overline{G}_+}.
\end{itemize}
\end{proposition}

{\bf Proof :}\quad
(i), (ii) Recall that for an arbitrary function $~f$,
\begin{equation}\label{iR}
\left(f^{\#}(s)\right)_{|_{s\in{i}\R}}=
\left(f(s)_{|_{s\in{i}\R}}\right)^*.
\end{equation}
Thus, whenever $~f\in{\mathcal Odd}~$ it maps $~i\R~$ to itself and
hence it is a $~{\mathcal G}{\mathcal P}~$ function. Next, from
\eqref{evenoddDef} and \eqref{iR} it follows that $~f_{\rm even}$,
the even part of an arbitrary $~f$, maps $~i\R~$ to $~\R$. Thus,
if $~f=f_{\rm even}+f_{\rm odd}~$ maps  $~i\R~$ to $~i\R$, it follows
that $~f_{\rm even}(s)_{|_{s\in{i}\R}}\equiv 0$, which in turn
means that $~f_{\rm even}(s)\equiv 0$.

Item (iii) follows from item (ii).

Item (iv) follows from item (i) together with \eqref{po} and
Theorem \ref{Th:gp}.\qed
\vskip 0.2cm

We now explore the structure of $~{\mathcal Odd}~$ from a
geometric point of view.

\begin{observation}\label{inner}
For a given set of rational functions $~G~$ ($G\not=\{ 0\}$),
denote by $~F~$ and $~H~$ the following sets,
$$\begin{matrix}
F:=\{g^{-1}~:~g\in G\}
&~&~&~&H:=\{~h=g^2~:~g\in G\}
\end{matrix}$$
$G\subset{\mathcal Odd}~$ if and only if
\begin{equation}\label{Inner}
{\rm Re}\left(\left(f(s)h(s)\right)_{|_{s\in i\R}}\right)\equiv 0
\quad\quad\quad\quad\forall f\in F\quad\forall h\in H.
\end{equation}
\end{observation}

{\bf Proof :}\quad 
The claim relies on Proposition \ref{gpo} item (ii).
First, if $g_a(s)$ and $g_b(s)$ are odd
functions then
$~{\rm Im}\left(g_a^2(s)_{|_{s\in{i\R}}}\right)\equiv 0$
and \mbox{${\rm Re}\left(g_b^{-1}(s)_{|_{s\in{i\R}}}
\right)\equiv 0$.} Thus, 
\mbox{${\rm Re}\left(\left(g_a^2(s)g_b^{-1}(s)
\right)_{|_{s\in{i\R}}}\right)\equiv 0$,} i.e.
\eqref{Inner} is satisfied.
\vskip 0.2cm

Conversely, if \eqref{Inner} is satisfied for all $~f\in F~$
and all $~h\in H$, it in particular holds for $f=g^{-1}~$ and
$h=g^2$ with the same $g$. This implies that $fh=g$, which
means $g\in{\mathcal Odd}$.
\qed\mbox{}
\vskip 0.2cm

Note that for a pair of functions $~f, h~$ one can define a
function-valued inner product
\mbox{$<f,h>~:={\rm Re}\left(f(s)h^{\#}(s)\right)$}
(in the sense that \mbox{$<rf,h>=r<f,h>$} for $~r\in\R)$. Now,
as $~\left(f(s)h(s)\right)_{|_{s\in i\R}}=\left(f(s)h^{\#}(s)
\right)_{|_{s\in i\R}}$, equation \eqref{Inner} can be
written as $~<f,h>_{|_{s\in i\R}}\equiv 0$. Namely, one can
interpret Observation \ref{inner} as saying that
$~0\not\equiv g\in{\mathcal Odd}~$ is equivalent to having
the restriction to the imaginary axis of $~g~$
and $~g^2$, orthogonal in the above inner product.
\vskip 0.2cm

Based on Observation \ref{ConvexPartitioning} and Proposition
\ref{gpo} we can now introduce a ~{\em convex} partitioning of all
$~{\mathcal Odd}~$ functions. For a fixed
$~g\in\overline{{\mathcal G}_+}$ the set $~\mathcal{GP}_g~$ was
defined in \eqref{eq:GP_g}. We now consider the set odd functions
in it:
$${\mathcal Odd}_g:=\mathcal{GP}_g\cap{\mathcal Odd}.$$
This subset is given by
\begin{equation}\label{Odd_g}
{\mathcal Odd}_{g}=\bigcup\limits_{a_o\geq 0, a_j>0, r_j\in\R}~g
\left( ir_o+a_os+\sum\limits_{j\geq 1}\frac{a_j}{s-ir_j}
\right)g^{\#}.
\end{equation}

\begin{observation}\label{ConvexPartitioningOdd}
Let $~{\mathcal Odd}_g~$ be as in \eqref{Odd_g}. Then,
$${\mathcal Odd}
=\bigcup\limits_{g\in\overline{{\mathcal G}_+}}
{\mathcal Odd}_g~.$$
\end{observation}

Recall that in \cite{YS} it was shown that (up to possibly
compromising the minimal degree interpolating function)
without loss of generality the classical Nevanlinna-Pick
interpolation may be confined to
$\mathcal{P}\cap{\mathcal Odd}$. Similarly, for a fixed
$~g\in\overline{{\mathcal G}_+}$ interpolation of 
$~\mathcal{GP}_g~$ functions may be confined to
${\mathcal Odd}_g~$ functions. This is illustrated next.

\begin{example}\label{Ex:Odd_g}
{\rm
In Example \ref{Ex:GenNevPickGP_g} we studied a variant of
Nevanlinna-Pick interpolation problem within the set
$~\mathcal{GP}_g~$ with $~g(s)=\frac{4}{7-3s}$, and looked
for a function $~\psi$ so that
\[
\psi(s)_{|_{s=1}}=1
\quad\quad{\rm and}\quad\quad
\psi(s)_{|_{s=2}}=4.
\]
It turned out that this was equivalent for a classical
Nevanlinna-Pick interpolation problem of searching
$p\in\mathcal{P}$ so that
\begin{center}
$p(s)_{|_{s=1}}=${\mbox{\tiny$\frac{5}{2}$}}$
\quad\quad{\rm and}\quad\quad
p(s)_{|_{s=2}}=${\mbox{\tiny$\frac{13}{4}$}}~.
\end{center}
Following the above analysis, looking for
$\psi\in{\mathcal Odd}_g$ is equivalent to restricting 
the classical Nevanlinna-Pick search for
$~p\in{\mathcal P}\cap{\mathcal Odd}$. From \cite{YS}
it follows that this restriction, does not limit the
solvability of the problem. In fact, there are still
infinitely many solutions.

We here mention two solutions $~p_a(s)=\frac{3s(s^2+9)} {4(s^2+2)}~$ (from Example
\ref{Ex:GenNevPickGP_g}) and \mbox{$p_d(s)=\frac{1}{6}(8s+\frac{7}{s})$.} 
The resulting interpolating functions in $~{\mathcal Odd}_g~$ are
\begin{center}
$\begin{matrix}
\psi_a(s)&=&g(s)p_a(s)g^{\#}(s)&=&
\frac{12s(s^2+9)}{(s^2+2)(49-9s^2)},\\
\psi_d(s)&=&g(s)p_d(s)g^{\#}(s)&=&
\frac{8(8s^2+7)}{3s(49-9s^2)}~.
\end{matrix}$
\end{center}
The function $\psi_a(s)$ is of degree four while $\psi_d(s)$ is of
degree three. 
\vskip 0.2cm

Following item (iii)(a) in Proposition \ref{gpo}, taking
$~r\in\R~$ as parameter, $~r\psi_a(s)+(1-r)\psi_d(s)~$
forms a variety of interpolating functions within the same
${\mathcal Odd}_g$.
\vskip 0.2cm

Finally, a straightforward use of \eqref{po} reveals that the
above $\psi_d$ is 
%
%
a minimal degree interpolation function in ${\mathcal Odd}_g$.
}\qed
\end{example}

\section{Even generalized positive functions}
\label{sec:even}\setcounter{equation}{0}

In Proposition \ref{gpo} we showed that all ~{\em odd}~ functions
are generalized positive functions and characterized them. We now
characterize \mbox{$\mathcal{GPE}:=\mathcal{GP}\cap{\mathcal Even}$,}
the subset of even functions within $~\mathcal{GP}$.
We shall use again the convention that $~\prod\limits_1^0=1$.

\begin{proposition}\label{gpeChar}
The following are equivalent

\begin{itemize}
\item[(i)~~~~]{}$\psi\in\mathcal{GPE}.$

\item[(ii)~~~]{}$\psi(s)=c\cdot\frac
{\prod\limits_{j=1}^m
\left(1-\alpha_j(1+(s-i\beta_j )^2)\right)}
{\prod\limits_{k=1}^n
\left(1-\gamma_k(1+(s-i\delta_k)^2)\right)}~$
with $~c>0$, $~\alpha_j, \gamma_k\in(0,~1]$,
$~\beta_j, \delta_k\in\R$.

\item[(iii)~~]{}$\psi(s)=g(s)g^{\#}(s)~$ for some
rational $~g(s)$.

\item[(iv)~~]{}
$\psi(s)~$ maps $~i\R~$ to $~\overline{\R_+}$.

\item[(v)~~~]{}$\psi\in{\mathcal G}{\mathcal P}~$ maps
$~i\R~$ to $~\R$.
\end{itemize}
\end{proposition}

\noindent
{\bf Proof}\quad
Any even function maps $~i\R~$ to $~\R$, thus
$~(i)~$ implies $~(v)$.

$(v)~\Longleftrightarrow~(iv)$. From $~(v)~$
it follows that $~\psi~$ maps $~i\R~$ to
$~\overline{\C_+}\cap\R$, thus in fact to
$~\overline{\R_+}$, so this part is established.

$(iv)~\Longleftrightarrow~(iii)$. From Theorem \ref{Th:gp}
together with with fact that\\
\mbox{$g(s)\psi(s)g^{\#}(s)_{|_{s\in{i\R}}}=
g(s)\psi(s)\left(g(s)\right)^*$,} it follows that
\mbox{$\psi(s)=g(s)p(s)g^{\#}(s)$,}
where $~p\in{\mathcal P}~$ maps $~i\R~$ to $~\overline{\R_+}$.
Now, up to non-negative scaling, $~p(s)\equiv 1~$ in whole $~\C$,
is the only function which achieves that.

$(iii)~\Longrightarrow~(ii)$. Denote,
$~g(s)=\tilde{c}\frac{\prod\limits_{j=1}^m(s-z_j)}
{\prod\limits_{k=1}^n(s-p_k)}~$ with $~\tilde{c}, z_j, p_k\in\C.$
Thus,
\[
\begin{split}
g(s)g(s)^{\#}&=
|\tilde{c}|^2\cdot\frac{\prod\limits_{j=1}^m\left(-s^2+s(z_j-z_j^*)
+|z_j|^2\right)}{\prod\limits_{k=1}^n\left(-s^2+s(p_k-
p_k^*)+|p_k|^2\right)}\\
&=|\tilde{c}|^2\cdot\frac{\prod\limits_{j=1}^m \left( ({\rm
Re}(z_j))^2-(s-i{\rm Im}(z_j))^2\right)}{
\prod\limits_{k=1}^n\left( ({\rm Re}(p_k))^2-(s- i{\rm
Im}(p_k))^2\right)}\\
&=|\tilde{c}|^2\cdot\frac{\prod\limits_{j=1}^m
\left(1+({\rm Re}(z_j))^2-(1+
(s-i{\rm Im}(z_j))^2)\right)}
{\prod\limits_{k=1}^n\left(1+({\rm Re}(p_k))^2-
(1+(s- i{\rm Im}(p_k))^2)\right)}\\
&=|\tilde{c}|^2\cdot\frac{\prod\limits_{j=1}^m
\left(1+({\rm Re}(z_j))^2\right)}
{\prod\limits_{k=1}^n\left(1+({\rm Re}(p_k))^2\right)}
\cdot
\frac{\prod\limits_{j=1}^m\left(1-
\frac{1}{1+({\rm Re}(z_j))^2}
(1+(s-i{\rm Im}(z_j))^2)\right)}
{\prod\limits_{k=1}^n\left(1-
\frac{1}{1+({\rm Re}(p_k))^2}(1+(s- i{\rm
Im}(p_k))^2)\right)}~,
\end{split}
\]
Denoting $~c:=|\tilde{c}|^2\cdot\frac{\prod\limits_{j=1}^m
\left(1+({\rm Re}(z_j))^2\right)}
{\prod\limits_{k=1}^n\left(1+({\rm Re}(p_k))^2\right)}~$,
$~\alpha_j=\frac{1}{1+({\rm Re}(z_j))^2}~$,
$\beta_j={\rm Im}(z_j)$,
$~\gamma_k=\frac{1}{1+({\rm Re}(p_k))^2}~$ and
\mbox{$\delta_k={\rm Im}(p_k)$}~
completes the construction.\\
As trivially $~(ii)~\Longrightarrow~(i)$, the claim is
established. \qed\mbox{}\\

\noindent
From Proposition \ref{gpeChar} it follows that $~\psi(s)~$ is in
$\mathcal{GPE}$ may be characterized as \mbox{$\psi=gg^{\#}$.}
Note however that this
factorization is non-unique, namely, one can have $g_1\not=g_2$
and still $g_1g_1^{\#}=g_2g_2^{\#}$. A characterization of all
these factorizations is given in \cite{F}.

One can now state several properties of
$\mathcal{GPE}$, the subset of even functions within $~\mathcal{GP}$.

\begin{proposition}\label{gpEv}
Let $~\psi~$ be a rational function. The following statements
are true.

\begin{itemize}
\item[(i)~~~~]{}$\psi\in\mathcal{GP}\Longleftrightarrow
\psi_{\rm even}\in\mathcal{GP}$.

\item[(ii)~~~]{} $\mathcal{GPE}~$
is a sub{\bf cic} of $~\mathcal{GP}$.

\item[(iii)~~]{}$\mathcal{GP}\cap{\mathcal Even}~$ is a
multiplicative group.

\item[(iv)~~]{}$g\in\mathcal{GPE}~~
\Longleftrightarrow~~g\cdot\mathcal{GP}
\subset\mathcal{GP}$.

\item[(v)~~~]{}For arbitrary $~\mathcal{GPE}~$ 
functions $~g_1(s)g_1^{\#}(s)~,~\ldots~,~g_m(s)g_m^{\#}(s)$,
there always exists $~\hat{g}\hat{g}^{\#}\in\mathcal{GPE}~$
so that
$$\sum\limits_{j=1}^mg_j(s)g_j^{\#}(s)=\hat{g}(s)\hat{g}^{\#}(s).$$
Moreover, one can always take $~\hat{g}\in\overline{{\mathcal G}_+}$.

\item[(vi)~~~]{}Let $~\psi~$ be the composition function
\mbox{$\psi(s):=p(g(s))$} where $g\in\mathcal{GPE}$ and
$~p\in{\mathcal P}$. Then, $\psi\in\mathcal{GP}$. If in
addition $~p~$ leaves the real axis invariant (e.g. $~p~$ is
real), then $~\psi\in\mathcal{GPE}$.
\end{itemize}
\end{proposition}

\noindent
{\bf Proof}\quad
(i)~ This follows from the fact that
$~\psi\in{\mathcal G}{\mathcal P}\Longleftrightarrow
\psi^{\#}\in{\mathcal G}{\mathcal P}.$\\
(ii) Recall that the set $~\mathcal{GP}~$ is a Convex
Invertible Cone, the claim is immediate from Proposition
\ref{gpeChar}(iv). Alternatively, $\mathcal{GPE}$ is a
non-empty intersection of two {\bf cic}s and thus a
sub{\bf cic}, \cite[Observation 2.1]{CL4}.\\
(iii) Is immediate from items (ii) or (iv) in
Proposition \ref{gpeChar}.\\
(iv) If $~f~$ is a generalized positive function, from item (i) we
know that \mbox{$f_{\rm even}\in\mathcal{GP}$.} Now if
$~g\in\mathcal{GPE}~$ from item (ii) it follows that
\mbox{$(f_{\rm even}g)\in\mathcal{GPE}$.} Next,
from Proposition \ref{EvenOdd}(ii) it follows that
$~(fg)_{\rm even}=f_{\rm even}g~$ and hence,
$~(fg)_{\rm even}\in\mathcal{GPE}$. Using again
item (i) implies, $~(fg)\in\mathcal{GP}\cap{\mathcal Even}$.

For the other direction all we need to show is that if $~g~$ is
a non-even function within $~\mathcal{GP}$, one can always find
$~f\in\mathcal{GP}~$ so that $~(fg)\not\in\mathcal{GP}$. Indeed
assume that $~g~$ is so that $~g(s)_{|_{s=i\omega_o}}=a+ib~$
with $~a\geq 0~$ and $~0\not=b\in\R~$ for some $~\omega_o\in\R$.
Then, taking the (odd) $~\mathcal{GP}~$ function
$~f=b(\omega_os+i)$ reveals that
$~{\rm Re}\left(f(s)g(s)_{|_{s=i\omega_o}}\right)=
-b^2(1+\omega_o^2)~$ and thus $~gf\not\in\mathcal{GP}$.\\
(v) Is immediate from item (ii) here together with item
(iii) in Proposition \ref{gpeChar}.\\
(vi) By construction $~g~$ maps $~i\R~$ to $~\overline{\R_+}~$
and in turn $~p~$ maps $~\overline{\R_+}~$ to $~\overline{\C_+}$,
thus, $~\psi~$ maps  $~i\R~$ to $~\overline{\C_+}$. If in
addition $~p~$ maps $~\overline{\R_+}~$ to $~\overline{\R_+}$,
$~\psi~$ maps  $~i\R~$ to $~\overline{\R_+}$, so the claim is
established.\qed\mbox{}\\

The convexity of the set $\mathcal{GP}$ and of its subset of
$\mathcal{GPE}$ functions, see item (v) in Proposition
\ref{gpEv}, may be exploited to introduce a straightforward
scheme of solving the Nevanlinna-Pick interpolation problem.

\begin{example}\label{ex:InterpGPE}
{\rm
For simplicity we consider the real case.

{\bf a}. We first look for a real polynomial $gg^{\#}\in\mathcal{GPE}$
so that
\begin{equation}\label{GPEinterp}
{g(s)g(s)^{\#}}_{|_{s=\pm 1}}=1,
\quad\quad\quad\quad
{g(s)g(s)^{\#}}_{|_{s=\pm 2}}=4,
\quad\quad\quad\quad
{g(s)g(s)^{\#}}_{|_{s=\pm 3}}=9.
\end{equation}
(Note that $f(s)=s^2$ is an interpolating polynomial, but
only $-f\in\mathcal{GPE}$).

Consider the following real $\mathcal{GPE}$ polynomials,
\begin{center}
$g_1(s)g_1^{\#}(s)=(4-s^2)(9-s^2)(${\mbox{\tiny$\frac{25}{24}$}}
$-s^2)$

$g_2(s)g_2^{\#}(s)=(1-s^2)(9-s^2)${\mbox{\tiny$\frac{1}{9}$}}
$(${\mbox{\tiny$\frac{8}{5}$}}$-s^2)$

$g_3(s)g_3^{\#}(s)=(1-s^2)(4-s^2)${\mbox{\tiny$\frac{1}{360}$}}
$(90-s^2)$
\end{center}
It is easy to verify that,
\[
\begin{matrix}
s&=&\pm 1&~&\pm 2&~&\pm 3\\
g_1g_1^{\#}&=&1&~&0&~&0\\
g_2g_2^{\#}&=&0&~&4&~&0\\
g_3g_3^{\#}&=&0&~&0&~&9
\end{matrix}
\]
Thus, using item (v) in Proposition \ref{gpEv},
$gg^{\#}=g_1g_1^{\#}+g_2g_2^{\#}+g_3g_3^{\#}$ is a real
$\mathcal{GPE}$ polynomial satisfying \eqref{GPEinterp}.

{\bf b}. Using part ~{\bf a}, we now look for a real polynomial
$\psi\in\mathcal{GP}$ so that
\begin{equation}\label{GPinterp}
\psi(s)_{|_{s=1}}=1,\quad\quad\quad\quad
\psi(s)_{|_{s=2}}=4, \quad\quad\quad\quad
\psi(s)_{|_{s=3}}=9,
\end{equation}
(with no constraints on $\psi(s)_{|_{s=-1,-2,-3}}$).
For $j=1,2,3$ we now construct real $\mathcal{GP}$ polynomials
of the form $\psi_j=g_jp_jg_j^{\#}$ with $g_jg_j^{\#}$ from
part ~{\bf a}~ and $p_j\in\mathcal{P}$ are of the form
$\frac{a_j+j}{a_j+s}$ with $a_j>0$ is so that one of the
roots of $g_jg_j^{\#}$ is canceled. Indeed take,
\begin{center}
$p_1(s)=${\mbox{\tiny$\frac{\frac{5}{2\sqrt{6}}+1}{\frac{5}
{2\sqrt{6}}+s}$}}$~,
\quad\quad\quad\quad
p_2(s)=${\mbox{\tiny$\frac{\frac{2\sqrt{2}}{\sqrt{5}}+2}
{\frac{2\sqrt{2}}{\sqrt{5}}+s}$}}$~,
\quad\quad\quad\quad
p_3(s)=${\mbox{\tiny$\frac{3\sqrt{10}+3}{3\sqrt{10}+s}$}}$~.$
\end{center}
Thus, one obtains,
\begin{center}
$\psi_1(s)=(4-s^2)(9-s^2)(${\mbox{\tiny$\frac{5}{2\sqrt{6}}$}}
$+1)(${\mbox{\tiny$\frac{5}{2\sqrt{6}}$}}$-s)$

$\psi_2(s)=(1-s^2)(9-s^2)
${\mbox{\tiny$\frac{2}{9}$}}$(${\mbox{\tiny$\frac{\sqrt{2}}
{\sqrt{5}}$}}$+1)(${\mbox{\tiny$\frac{2\sqrt{2}}{\sqrt{5}}$}}
$-s)$

$\psi_3(s)=(1-s^2)(4-s^2)${\mbox{\tiny$\frac{1}{120}$}}
$(\sqrt{10}+1)(3\sqrt{10}-s)$
\end{center}
It is easy to verify that,
\[
\begin{matrix}
s&=&1&~&2&~&3\\
\psi_1&=&1&~&0&~&0\\
\psi_2&=&0&~&4&~&0\\
\psi_3&=&0&~&0&~&9.
\end{matrix}
\]
Thus, taking $\psi=\psi_1+\psi_2+\psi_3$ satisfies \eqref{GPinterp}.
\vskip 0.2cm

Roughly speaking, the simplicity of this scheme of constructing
interpolating functions, comes on the expense of high degree.
}\qed\end{example}

As already mentioned, the fact that the set ${\mathcal P}$ is
a convex invertible cone {\bf cic}, is classical, see e.g.
\cite[5.6]{Be} and item (i) in Proposition \ref{MaxG}. Recall
that real $~{\mathcal P}~$ functions are identified with the
driving point impedance of ~\mbox{R-L-C}~ electrical circuits
\cite{AV}, \cite{Be}, \cite{brune}, \cite{Br}. In fact, in the
framework of \mbox{R-L-C} electrical circuits the
three {\bf cic} operations of positive scaling, summation and
inversion have the physical interpretation of transformer ratio,
series connection of impedances and impedance/admittance duality,
respectively. Moreover, recall that
$~{\mathcal P}\cap{\mathcal Even}~$ is associated with resistive
circuits and $~\mathcal{PO}~$ with reactive
networks (L-C). However, not every network can be realized as a
series connection of a resistive and a reactive circuits. Namely,
it is only over $~\mathcal{GP}~$ that the partitioning of a
positive function into even and odd parts is always possible.
\vskip 0.2cm

Recall that in contrast to \eqref{po}, the set
$~{\mathcal P}\cap{\mathcal Even}~$ is almost empty, i.e. up to
positive scaling it consists of a single function, $~p(s)\equiv 1$.
We now show that if one is interested in even-odd partitioning of
functions, the set $~\mathcal{GP}~$ is closed, while
its subset of positive functions is not. Namely, $~p_{even}$, the
even part of a positive function $~p$, is either a non-negative
constant or not a positive function. One can only guarantee that
$~p_{\rm even}\in\mathcal{GPE}$. This is illustrated next.

\begin{example}
{\rm
Consider the positive (real) function: $~\psi(s)=\frac{1}{1+s}~$
defined in the whole $\C$. Then $~p_{\rm even}(s)=\frac{1}{1-s^2}~$
and $~p_{\rm odd}(s)=\frac{-2s}{1-s^2}~$ are so that
$p_{\rm even}\in\mathcal{GPE}~$ and
$~p_{\rm odd}\in\mathcal{GP}\cap{\mathcal Odd}$,
but
neither $~p_{\rm even}~$ nor $~p_{\rm odd}~$ are positive.}
\qed\end{example}

We conclude this section by introducing yet another factorization
of $~\mathcal{GP}~$ functions through odd functions.

\begin{observation}
$\psi\in\mathcal{GP}~$ if and only if there exist
$~f,g\in{\mathcal Odd}~$ so that $~\psi_{\rm even}=-f^2~$ and
$~\psi_{\rm odd}=g$.
\end{observation}

{\bf Proof :}\quad
Since $~\psi_{\rm odd}\in{\mathcal Odd}$, we only need to
show that $~\psi\in{\mathcal G}{\mathcal P}~$  if and only if
\mbox{$\psi_{\rm even}=-f^2$} for some $~f\in{\mathcal Odd}$.
Now from item (i) in Proposition \ref{gpEv} it follows that
having $~\psi\in\mathcal{GP}~$ is equivalent to
$~\psi_{\rm even}\in\mathcal{GP}$. From item (iv) in
Proposition \ref{gpeChar} this in turn is equivalent to
$~\psi_{\rm even}~$ mapping $~i\R~$ to $~\overline{\R_+}$.
Using item (ii) from Proposition \ref{gpo} completes the
proof.
\qed

\section{Generalized bounded functions}\label{sec:bounded}
\setcounter{equation}{0}

Recall that a function $~f_b(s)~$ is called ~{\em bounded},
denoted by $~f_b\in{\mathcal B}$, (commonly the real case is
addressed) if it analytically maps $~\C_+~$ to the closed unit
disk, see e.g. \cite[Chapter 7]{AV}, \cite[Section 6.5]{Be}
and $~f_{gb}(s)~$ is {\em generalized bounded}~
$~f_{gb}\in\mathcal{GB}~$ if it maps $~i\R~$ to the closed
unit disk, see e.g. \cite{DDGK}. It is well known that 
through the Cayley transform one can identify positive
functions with bounded functions, namely
\begin{equation}\label{cayley1}
\begin{matrix}
f_b(s)=\dfrac{1-p(s)}{1+p(s)}&~&p\in\mathcal{P},&~&~&~&
f_{gb}(s)=\dfrac{1-\psi(s)}{1+\psi(s)}&~&
\psi\in\mathcal{GP}. 
\end{matrix}
\end{equation}
Nevertheless, we here focus on the less obvious analogies.
In Proposition \ref{gb1} and Corollary \ref{gb2} below we introduce
two representations of $\mathcal{GB}$ functions.

\begin{proposition}\label{gb1}
A rational $f_{gb}(s)$ is a generalized bounded function
if and only if it is of the form $f_{gb}(s)=f_b(s)/\beta(s)$,
where $f_b\in\mathcal B$ and where $\beta(s)$ is a finite
Blaschke product.
\end{proposition}

{\bf Proof:} We first note that since $f_{gb}$ is bounded
on the imaginary axis, all its singularities there are
removable. Let $w_1,\ldots, w_\ell$ be the poles of
$f_{gb}$ in ${\mathbb C}_+$, and consider the function
\[
f_b(s)=f_{gb}(s)\prod_{j=1}^\ell \frac{s-w_j}{s+w_j^*}~.
\]
The function $f_{gb}$ is analytic and bounded by $1$ in
modulus in ${\mathbb C}_+$, as is seen for example by the
maximum modulus principle, or by direct inspection. We thus
have the result with
\[
\beta(s)=\prod_{j=1}^\ell \frac{s-w_j}{s+w_j^*}~.
\]
The converse is clear.
\mbox{}\qed\mbox{}\\

One can characterize generalized bounded functions through
the associated kernel.

\begin{corollary}
A rational $~f(s)~$ is a generalized bounded function if and
only if the kernel
\[
k_f(s,w)=\dfrac{1-f(s)f(w)^*}{s+w^*}
\]
has a finite number of negative squares in ${\mathbb
C}_+\smallsetminus \left\{w_1,\ldots, ,w_\ell\right\}$, where the
$w_j$ denote the poles of $f$ in $\mathbb C_+$.
\end{corollary}

{\bf Proof:} Assume that $f=f_{gb}=f_{b}/\beta$. As proved in a
more general context in \cite[Theorem 6.6, p. 132]{ad3}, one
direction follows from the equality
\[
k_{f_{gb}}(s,w)=\frac{1}{\beta(s)}
\left\{k_{f_{b}}(s,w)-k_{\beta}(s,w)\right\}\frac{1}{\beta(w)^*}~,
\]
see for instance the formula on top of page 134 in \cite{ad3}.\\

The converse is just a particular case of the above mentioned
result of Kre\u{\i}n Langer \cite[Theorem 3.2]{KL1}. A direct
proof for the rational case can also be given, but will be
omitted here.
\mbox{}\qed\mbox{}\\

We now turn to another representation of
$f_{gb}\in{\mathcal{GB}}$. Obviously, (generalized) bounded
functions and (generalized) positive functions are related
through the Cayley transform \eqref{cayley1}. We now introduce
an adapted version of this characterization. To this end,
recall (Proposition \ref{gpEv}) that the set
$~\mathcal{GP}\cap{\mathcal Even}~$ is characterized by functions
of the form $~gg^{\#}$. From the above discussion one has the
following.

\begin{corollary}\label{gb2}
A rational function $~f_{gb}(s)~$ is generalized bounded if
and only if it admits a representation,
\begin{equation}\label{gb3}
f_{gb}(s)=\left(g(s)g^{\#}(s)-p(s)\right)\left(
g(s)g^{\#}(s)+p(s)\right)^{-1}~,
\end{equation}
for some $~p\in{\mathcal P}~$
and some $~g(s)\in{\overline{G}_+}$, \eqref{overline{G}_+}
\end{corollary}

\noindent
{\bf Proof}\quad Indeed, from \eqref{char} and
\eqref{cayley1} it follows that $~f_{gb}\in\mathcal{GB}$
can be written as,
$$\begin{matrix}
f_{gb}&=&\left(1-\psi\right)\left(1+\psi\right)^{-1}=
\left(1-gpg^{\#}\right)\left(1+gpg^{\#}\right)^{-1}\\~&
=&\left(g\left( (g^{\#}g)^{-1}-p
\right)g^{\#}\right)
\left(g\left(
(g^{\#}g)^{-1}+p\right)g^{\#}\right)^{-1}\\~&=&
\left( (g^{\#}g)^{-1}-p\right)
\left( (g^{\#}g)^{-1}+p\right)^{-1}.
\end{matrix}$$
Now, from Proposition \ref{gpeChar} and Proposition \ref{gpEv}
it follows that $~\psi\in\mathcal{GP}\cap{\mathcal Even}~$
is equivalent to $~\psi(s)=\left(g^{\#}(s)g(s)\right)^{-1}$,
so up to inversion,
the claim is established. \qed\mbox{}\\

It is interesting to compare Corollary \ref{gb2}
with $~f_{gb}~$ in \eqref{cayley1}
\vskip 0.2cm

We conclude by pointing out that there is a structural
difference between $\mathcal{GP}$ and $\mathcal{GB}$
functions. One may be tempted to try to mimic, in the
framework of $\mathcal{GB}$ functions, the convex
partitioning of sets of functions of the form
$\mathcal{GB}_g$ in the spirit of \eqref{eq:GP_g} and
Observations \ref{ConvexPartitioning} and
\ref{ConvexPartitioningOdd}, sets where in $~\C_+~$ the
poles and zeroes are fixed.
However, unfortunately this is no longer true. This also
prevents us from mimicking the interpolation over
$\mathcal{GP}_g$ to $\mathcal{GB}_g$ functions.

Indeed, for a given $~g~$ in \eqref{gb3} define the set
$\mathcal{GB}$ function,
\begin{equation}\label{eq:GB_g}
\mathcal{GB}_g:=
\left\{\left(g(s)g^{\#}(s)-p(s)\right)\left(
g(s)g^{\#}(s)+p(s)\right)^{-1}~:~p\in{\mathcal P}~
\right\}.
\end{equation}
In the following example we show that in contrast to the set
$\mathcal{GP}_g$ in \eqref{eq:GP_g}, within the set
$\mathcal{GB}_g$ neither the poles nor the zeroes in $~\C_+~$
are fixed.

\begin{example}
{\rm
Fix in \eqref{eq:GB_g} $~g\in{\mathcal G}_+~$, see
\eqref{eq:G_+}, namely, \mbox{$g(s)=\frac{n}{d}~$} with
$~n=\prod\limits_{j=1}^l(s-z_j)~$ and
$~d=\prod\limits_{k=1}^q(s-\pi_k)$, where
$~z_1,~\ldots~,~z_l~$ and $~\pi_1,~\ldots~,~\pi_q~$
are prescribed (not necessarily distinct) points in $~\C_+$.
Let now, \mbox{$\tilde{n}=\prod\limits_{j=1}^l(s-z_j-\epsilon_j)~$}
and $~\tilde{d}=\prod\limits_{k=1}^q(s-\pi_k-\delta_k)~$
with $~\epsilon_j\geq 0$, $~\delta_k\geq 0$,
\mbox{$1>>\sum\limits_{j=1}^l\epsilon_j>0$,}
$~1>>\sum\limits_{k=1}^q\delta_k>0$. Let $~p_1:=\frac{n}{\tilde{n}}~$
and $~p_2:=\frac{\tilde{d}}{d}$. By construction, each of the
functions $~p_1, p_2~$ analytically maps $~\C_+~$ to a neighborhood of
the point +1 and hence both are positive. Thus both
$$\begin{matrix}
f_{gb,1}&:=&\left(g(s)g^{\#}(s)-p_1(s)\right)\left(
g(s)g^{\#}(s)+p_1(s)\right)^{-1}\\
f_{gb,2}&:=&\left(g(s)g^{\#}(s)-p_2(s)\right)\left(
g(s)g^{\#}(s)+p_2(s)\right)^{-1}\end{matrix}$$
are in $~\mathcal{GB}_g$.
However, in $~\C_+~$ they do not share the poles nor the zeroes.
Indeed, substitution yields $~f_{gb,1}=(\tilde{n}g^{\#}-d)(
\tilde{n}g^{\#}+d)^{-1}~$ and
$~f_{gb,2}=(ng^{\#}-\tilde{d})(ng^{\#}+\tilde{d})^{-1}$.
}
\qed
\end{example}

It should be pointed out that the above discussion reflects a
property of $\mathcal{GB}_g$ functions, independent of the
choice of the representation. Indeed, if the set in \eqref{eq:GB_g}
is substituted by the analogous one, based on Proposition \ref{gb1},
a conclusion, similar to that of the above example, is reached.
\vskip 0.2cm

It should be emphasized that rational generalized Nevanlinna
functions, mapping the real axis to the upper half plane, admit
a partitioning along the lines of Section \ref{sec:GP_g}. In
contrast, rational generalized Schur functions mapping the unit
circle to the unit disk, share the same difficulty as $\mathcal{GB}$
functions. This suggests that the complicated known scheme for
solving Nevanlinna- Pick interpolation problem for generalized
Schur functions, see e.g. \cite{Ball}, \cite{BH} and \cite{DD}
and for the single point with derivatives version, \cite{AADLW},
can not be simplified along the lines suggested in Example
\ref{Ex:GenNevPickGP_g}

\section{Future research}\label{sec:future}

In this work, part of ongoing research on $\mathcal{GP}$ functions,
we concentrated on exploring structural properties this set.
This opens the door for studying
various questions and we here mention sample of those.
First, in the framework of Nevanlinna-Pick interpolation
problem of scalar rational $\mathcal{GP}$ functions.

\begin{itemize}

\item{}Explore the question of minimal degree interpolating
functions.

\item{}We conjecture that if the Pick matrix $\Pi$ has $m$
negative eigenvalues then there exists an interpolating
function within a set $\mathcal{GP}_g$ where $g$ has $m$ poles
or $m$ zeroes in $\C_+$.

\item{}Parameterize all $\mathcal{GPE}$ interpolating functions.

\item{}Characterize the Nevanlinna-Pick interpolation
problems solvable by ${\mathcal Odd}$ and by $\mathcal{GPE}$ 
functions.
\end{itemize}

One can then look for generalizations. For example,

\begin{itemize}
\item{}Formalize the extension of the study of $\mathcal{GP}_g$
functions, to the cases of: (i) not necessarily rational,
(ii) matrix valued.

\item{}Formalize the extension of the study of the even-odd
partitioning of $\mathcal{GP}$ functions, to the cases of: (i) not
necessarily rational, (ii) matrix valued.
\end{itemize}

\end{document}